\title{Coupling and Bernoullicity \\ in random-cluster and Potts models}
\author{Olle H\"{a}ggstr\"{o}m\thanks{Mathematical Statistics, G\"oteborg 
University, 412 96 G\"oteborg, Sweden,
\texttt{olleh@math.chalmers.se}, 
\texttt{http://www.math.chalmers.se/\~{ }olleh/}}\
\thanks{Research partially supported by the Swedish Natural Science 
Research Council.}
\and Johan Jonasson\thanks{Department of Mathematics, Chalmers 
University of Technology, 412 96 G\"oteborg, Sweden,
\texttt{jonasson@math.chalmers.se}, 
\texttt{http://www.math.chalmers.se/\~{ }jonasson/}}
\and Russell Lyons\thanks{Department of Mathematics,  
Indiana University, Bloomington, IN 47405-5701, USA,
\texttt{rdlyons@indiana.edu}, 
\texttt{http://php.indiana.edu/\~{ }rdlyons/} and 
School of Mathematics, Georgia Institute of Technology, Atlanta, GA
30332-0160, USA,
\texttt{rdlyons@math.gatech.edu}, 
\texttt{http://www.math.gatech.edu/\~{ }rdlyons/}
}\ 
\thanks{Research partially supported by NSF grant DMS-9802663 and by
Microsoft Corp.}}

\documentclass[11pt]{article}
\usepackage{amsmath,amssymb,amsthm,amsfonts,latexsym}
\newtheorem {thm}{Theorem}[section]
\newtheorem {prop}[thm]{Proposition}
\newtheorem {lem}[thm]{Lemma}

\newtheorem {conj}[thm]{Conjecture}
\newtheorem {quest}[thm]{Question}

\def\Cox{\hfill \Box}
\def\Z{{\mathbb Z}}
\def\R{{\bf R}}
\def\P{{\bf P}}
\def\T{{\bf T}}
\def\Pt{{\sf Pt}}
\def\FPt{{\sf FPt}}
\def\WPt{{\sf WPt}}
\def\RC{{\sf RC}}
\def\FRC{{\sf FRC}}
\def\WRC{{\sf WRC}}
\def\C{{\cal C}}
\def\leqd{\stackrel{\cal D}{\preccurlyeq}}
\def\geqd{\stackrel{\cal D}{\succcurlyeq}}

\def\Aut{{\rm Aut}}
\def\f{{\rm free}}
\def\w{{\rm wired}}
\def\eps{\varepsilon}

\def\ed#1{[#1]}  
\def\kk(#1){\|#1\|}  
\def\kkk(#1){\|#1\|^*}  
\def\I#1{{\bf 1}_{\{#1\}}}  
\def\st{\, \colon \;}  

\def\Seq#1{\langle #1 \rangle}

\def\dom{{\cal D}}  

\def\thmenv#1#2#3{\begin{#1} \label{#1:#2} #3 \end{#1}}

\def\procl#1.#2 #3\endprocl{%
       \ifx#1t\thmenv{thm}{#2}{#3}\fi
       \ifx#1l\thmenv{lem}{#2}{#3}\fi
       \ifx#1p\thmenv{prop}{#2}{#3}\fi
       \ifx#1c\thmenv{cor}{#2}{#3}\fi
       \ifx#1d\thmenv{defn}{#2}{#3}\fi
       \ifx#1g\thmenv{conj}{#2}{#3}\fi
       \ifx#1q\thmenv{quest}{#2}{#3}\fi
       \ifx#1r\thmenv{rmk}{#2}{{\rm #3}}\fi
    }%

\def\rref#1.#2/{%
      \ifx #1sSection~\ref{sect:#2}\fi
      \ifx #1tTheorem~\ref{thm:#2}\fi  
      \ifx #1lLemma~\ref{lem:#2}\fi 
      \ifx #1cCorollary~\ref{cor:#2}\fi 
      \ifx #1pProposition~\ref{prop:#2}\fi 
      \ifx #1dDefinition~\ref{defn:#2}\fi
      \ifx #1gConjecture~\ref{conj:#2}\fi 
      \ifx #1qQuestion~\ref{quest:#2}\fi 
      \ifx #1e(\ref{eq:#2})\fi
      \ifx #1b(\cite{#2})\fi
        }

\def\rlabel #1 #2{\begin{equation} \label{eq:#1} #2 \end{equation}}

\def\proof{\medbreak\noindent{\it Proof.\enspace}}

\def\proofof #1.#2 {\medbreak\noindent
     {\it Proof of \rref #1.#2/.}\enspace}

\def\Qed{$\Cox$\medbreak}

\begin{document}

\maketitle

\begin{abstract}
An explicit coupling construction of random-cluster
measures is presented. 
As one of the applications of the construction, 
the Potts model on amenable Cayley graphs is shown to 
exhibit at every temperature the mixing property known as Bernoullicity. 
\end{abstract}

\section{Introduction} \label{sect:intro}

In the (ferromagnetic) {\bf Potts model}, spins (or colors) from the 
set $\{1, \ldots, q\}$ are assigned to the vertices of a graph $G=(V,E)$
randomly, in a way that favors configurations where 
many pairs of neighboring
vertices take the same spin value. More precisely, a spin configuration
$\xi \in \{1, \ldots, q\}^V$ is assigned probability proportional to
\[
\exp \Big( -2\beta \sum_{\ed{ x,y} \in E}
\I{\xi(x) \neq \xi(y)} \Big)
\]
where $\beta\geq 0$ is referred to as the inverse temperature parameter. The
case $q=2$ is known as the Ising model. 

The Potts model has received a considerable amount of attention in the
statistical mechanics and probability literature for several decades. In the
last decade, perhaps the most important tool for analyzing the Potts model 
has been the {\bf random-cluster model}, which is a kind of edge representation
of the Potts model. It was introduced by Fortuin and Kasteleyn \cite{FK}, and
has been heavily exploited in the study of Potts models since the seminal 
papers by Swendsen and Wang \cite{SW}, Edwards and Sokal \cite{ES}, and
Aizenman, Chayes, Chayes and Newman \cite{ACCN}. One of the main points of
working with the random-cluster representation, rather than directly
with the Potts model, is that questions about spin correlations in the latter
turn into questions about connectivity probabilities in the former, thereby
allowing powerful percolation techniques to come into play. Another 
interesting aspect of the random-cluster representation is that it makes sense
also for noninteger $q$. 

This paper is a contribution to the study of random-cluster and Potts models 
on infinite lattices. After recalling some necessary prerequisites in
Section \ref{sect:preliminaries}, we come in Sections \ref{sect:dynamics} 
and \ref{sect:Bernoullicity} to the two main purposes of this paper, which
are the following:
\begin{itemize}
\item
In Section \ref{sect:dynamics}, we present a useful device for the
analysis of random-cluster and Potts models, namely an explicit pointwise
{\bf dynamical construction} of random-cluster measures. The construction
provides natural couplings between random-cluster measures with
different parameter values or different boundary conditions. 
To some extent, this construction can be viewed as known and our presentation
of it can to the same extent be viewed as expository; it consists of 
putting together a few well-known ingredients from
Grimmett \cite{Gr}, Propp and Wilson \cite{PW}, and H\"aggstr\"om, Schonmann
and Steif \cite{HSS}. 
\item
In Section \ref{sect:Bernoullicity}, we apply the dynamical construction 
from the preceding section to show that the Potts model with fixed-spin
boundary condition on 
$\Z^d$ (and more generally on amenable Cayley graphs) exhibits a rather
strong mixing condition known as {\bf Bernoullicity}. Our proof
appears to be the simplest to date even in cases where the result was
known previously. 
\end{itemize}
Finally, some additional consequences of, and questions on, the 
dynamical construction are discussed in Section \ref{sect:remarks}. 

\section{Preliminaries} \label{sect:preliminaries}

The following subsections are devoted to recalling known material that
will be used in later sections. The random-cluster and Potts models are
introduced in Sections \ref{subsect:random-cluster} and \ref{subsect:Potts},
respectively. Before that, however, we recall some graph terminology in
Section \ref{subsect:graphs} and some basics on stochastic domination in
Section \ref{subsect:stoch_dom}. A general reference for this 
background material is Georgii, H\"aggstr\"om and Maes \cite{GHM}.

\subsection{Some graph terminology}  \label{subsect:graphs}

Let $G=(V,E)$ be a graph with vertex set $V$ and edge set $E$. We shall
always assume either that the graph is finite, or that it is countably
infinite and locally finite. 
An edge $e\in E$ will often be denoted $\ed{ x,y }$. 
The number of edges incident to a vertex $x$ is called the {\bf degree} of
$x$. For $W \subset V$, we define the (inner) {\bf boundary} $\partial W$ of
$W$ as 
\begin{equation} \label{eq:inner_boundary}
\partial W := \{x \in W: \, \exists y \in V \setminus W \mbox{ such that }
\ed{ x,y} \in E\} \, . 
\end{equation}
A {\bf graph automorphism} of $G$ is a bijective mapping 
$\gamma: \, V \rightarrow V$ with the property that for all $x,y \in V$,
we have $\ed{ \gamma x, \gamma y } \in E$ if and only if
$\ed{ x,y } \in E$. Write $\Aut(G)$ for the group of all
graph automorphisms of $G$. To each $\gamma\in \Aut(G)$, there is a
corresponding mapping $\tilde{\gamma}: E\rightarrow E$ defined
by $\tilde{\gamma}\ed{ x,y }:= \ed{ \gamma x, \gamma y }$.
The graph $G$ is said to be {\bf transitive}
if and only if for some (any) $x \in V$, one has that
for any $y \in V$ there exists $\gamma\in\Aut(G)$
such that $\gamma x = y$. 
One says that $G$ is {\bf quasi-transitive} if and only if for some 
finite subset $\{x_1,\ldots,x_n\}$ of $V$, one has that for any
$y \in V$ there exists $\gamma\in\Aut(G)$ such that 
$\gamma x_i = y$ for some $x_i$.

A probability
measure $\mu$ on $\{0,1\}^E$ is said to be {\bf automorphism invariant} if
for
any $n$, any $e_1, \ldots, e_n \in E$, any $i_1, \ldots, i_n \in \{0,1\}$,
and any $\gamma\in \Aut(G)$ we have
\begin{align*}
\mu\bigl(\{ X \in \{ 0, 1 \}^E &\st X(e_1)=i_1, \ldots, X(e_n)=i_n \}\bigr)
\\&= 
\mu\bigl(\{X \in \{ 0, 1 \}^E \st X(\tilde{\gamma}(e_1))=i_1, \ldots,
X(\tilde{\gamma}(e_n))=i_n \}\bigr) \, .
\end{align*}
In the sequel, we shall simplify the notation and omit the ``$ \{ X \in \{
0, 1 \}^E \st \} $" as used in the preceding equation.

A graph property that turns out to be important in many situations is
amenability: An infinite graph $G$ is said to be {\bf amenable} if
\[
\inf \frac{|\partial W|}{|W|}=0 
\,,
\]
where the infimum ranges over all finite $W\subset V$, and $| \cdot |$
denotes cardinality. There are
various alternative definitions of amenability of a graph that coincide
for transitive graphs (and more generally for graphs of bounded degree), but
not in general. 

For any graph $G$ and $x\in V$, define the {\bf stabilizer} $S(x)$ as the
set
of graph automorphisms that fix $x$, i.e.,
\[
S(x) := \{\gamma \in \Aut(G)  \st \gamma x = x \} \, .
\]
For $x,y \in V$, define
\[
S(x)y := \{ z \in V  \st \exists \gamma \in S(x) \mbox{ such that }
\gamma y = z \} \, .
\]
When $\Aut(G)$ is given the weak topology generated by its action on $V$, all
stabilizers are compact subgroups of $\Aut(G)$
because $G$ is locally finite and connected.
A transitive graph $G$ is said to be
{\bf unimodular} if for all $x, y\in V$ we
have the symmetry
\[
|S(x)y| = |S(y)x| \, .
\]

Another important class of graphs is the class of
Cayley graphs. If $\Gamma$ is a finitely generated group with 
generating set
$\{g_1, \ldots, g_n\}$, then the {\bf Cayley graph} associated
with $\Gamma$ and that particular set of generators is the (unoriented)
graph
$G=(V,E)$ with vertex set $V:=\Gamma$, and edge set
\[
E:= \{ \ed{ x,y } \st x,y \in \Gamma, \exists i \in 
\{1, \ldots, n\} \mbox{ such that } x g_i = y \} \, .
\]
Obviously, a Cayley graph is transitive, and 
furthermore it is not hard to show that it is unimodular.  
Most graphs that have been studied in percolation theory
are Cayley graphs. Examples include $\Z^d$ (which, with a slight abuse of 
notation, is short for the graph with vertex set $\Z^d$ and edges connecting
pairs of vertices at Euclidean distance $1$ from each other), and the
regular 
tree $\T_n$ in which every vertex has exactly $n+1$ neighbors. The graph
$\Z^d$
is amenable, while $\T_n$ is nonamenable for $n\geq 2$. Also studied are
certain nonamenable tilings of the hyperbolic plane 
(see, e.g., \cite{BS:hp} and \cite{HJL}), and further examples can be
obtained, e.g., by taking Cartesian products of two or more Cayley graphs. 

\subsection{Stochastic domination}  \label{subsect:stoch_dom}

Let $E$ be any finite or countably infinite set. (In our applications,
$E$ will be an edge set; hence the notation.) For two
configurations $\xi, \xi' \in \{0,1\}^E$, we write $\xi \preccurlyeq \xi'$
if
$\xi(e) \leq \xi'(e)$ for all $e\in E$. A function 
$f:\{0,1\}^E \rightarrow \R$ is said to be increasing if $f(\xi) \leq
f(\eta)$
whenever $\xi\preccurlyeq \eta$. For two probability measures $\mu$ and
$\mu'$
on $\{0,1\}$, we say that $\mu$ is {\bf stochastically dominated} by $\mu'$,
writing $\mu \leqd \mu'$, if
\begin{equation} \label{eq:def_stoch_dom}
\int_{\{0,1\}^E} fd\mu \leq \int_{\{0,1\}^E} fd\mu' 
\end{equation}
for all bounded increasing $f$. 

By a {\bf coupling} of $\mu$ and $\mu'$, or of two random objects $X$
and $X'$ with distributions $\mu$ and $\mu'$, we simply mean
a joint construction of two random objects with the prescribed distributions
on a common probability space. 

By Strassen's Theorem (see, e.g., \cite{Li}), $\mu \leqd \mu'$
is equivalent to the existence of a coupling $\P$ of two random
objects $X$ and $X'$ with distributions $\mu$ and $\mu'$,
such that $\P(X\preccurlyeq X')=1$. We call such a coupling a {\bf witness}
to the stochastic domination (\ref{eq:def_stoch_dom}).  

A useful tool for establishing stochastic domination is the well-known
Holley's Inequality. For $E' \subset E$ and $\xi \in \{0,1\}^E$, we
let $\xi(E')$ denote the restriction of $\xi$ to $E'$. 
\begin{lem}[Holley's Inequality]  \label{lem:Holley}
Let $E$ be finite, and let $\mu$ and $\mu'$ be probability measures
on $\{0,1\}^E$ that assign positive probability to all elements of
$\{0,1\}^E$.
Suppose that $\mu$ and $\mu'$ satisfy
\[
\mu(X(e)=1 \, | \,  X(E \setminus \{e\}) = \xi) \leq
\mu'(X(e)=1 \, | \, X(E \setminus \{e\}) = \xi')
\]
for all $e \in E$, and all $\xi, \xi' \in \{0,1\}^{E\setminus \{e\}}$ 
such that $\xi \preccurlyeq \xi'$. Then $\mu \leqd \mu'$. 
\end{lem}
This is not the most general form of Holley's Inequality, but one that
is sufficient for our purposes. 
For a proof, see, e.g., \cite{GHM} (Theorem 4.8).

We shall also need the notion of weak convergence of probability measures
on $\{0,1\}^E$, when $E$ is countably infinite. 
For such probability measures $\mu_1, \mu_2, \ldots$ 
and $\mu$, we say that $\mu$ is the (weak) limit of $\mu_i$ as 
$i \rightarrow \infty$ if
$\lim_{i\rightarrow\infty}\mu_i(A) = \mu(A)$ for all cylinder events $A$. 

\subsection{The random-cluster model}  \label{subsect:random-cluster}

Let $G=(V,E)$ be a finite graph. 
An element $\xi$ of $\{0,1\}^E$ will
be identified with the subgraph of $G$ that has vertex set $V$ and edge
set $\{e\in E \st \xi(e)=1\}$. An edge $e$ with $\xi(e)=1$
(resp.\ $\xi(e)=0$) is said to be open (resp.\ closed). 
A central quantity to the random-cluster model is the
number of connected components of $\xi$, which will be denoted $\kk(\xi)$.
We emphasize that in the definition of $\kk(\xi)$, 
isolated vertices in $\xi$ also count as connected components. 

The {\bf random-cluster measure} $\RC:=\RC^G_{p,q}$ (sub- and superscripts
will be dropped whenever possible) with parameters $p \in [0,1]$ and $q>0$, 
is defined as the probability measure on $\{0,1\}^E$ that to each
$\xi\in\{0,1\}^E$ assigns probability
\begin{equation} \label{eq:RC_def}
\RC(\xi):= 
\frac{q^{\kk(\xi)}}{Z} \prod_{e\in E} p^{\xi(e)}(1-p)^{1-\xi(e)} \, ,
\end{equation}
where $Z:=Z^G_{p,q}:= \sum_{\xi\in\{0,1\}^E} q^{\kk(\xi)}
\prod_{e\in E} p^{\xi(e)}(1-p)^{1-\xi(e)}$ is a normalizing constant making
$\RC$ a probability measure. 

When $q=1$, we see that all edges are independently open and closed with
respective probabilities $p$ and $1-p$, so that we get the usual 
i.i.d.\ bond percolation model on $G$. All other choices of $q$ yield
dependence between the edges. Throughout the paper, we shall assume (as
in most studies of the random-cluster model)
that $q\geq 1$. The main reason for doing so is that when $q\geq 1$, 
the conditional probability in eq.\ (\ref{eq:single_edge_cond_prob}) below
becomes increasing not only in $p$ but also in $\xi$, and this allows
some very powerful stochastic domination arguments, based on Holley's
Inequality (Lemma \ref{lem:Holley}), to come into play; these are not
available for $q< 1$. Furthermore, it is only random-cluster measures with
$q\in\{2,3,\ldots\}$ that have proved to be useful in the analysis of 
Potts models. 

It is immediate from the definition 
that if $X$ is a $\{0,1\}^E$-valued random object with
distribution $\RC$, then we have, for each $e = \ed{x,y} \in E$ and each 
$\xi\in\{0,1\}^{E\setminus\{e\}}$, that
\begin{equation} \label{eq:single_edge_cond_prob}
\RC\bigl(X(e)=1 \bigm| X(E\setminus\{e\})= \xi\bigr) =
\left\{
\begin{array}{ll}
p & \mbox{if } x\leftrightarrow y, \\
\frac{p}{p+(1-p)q} & \mbox{otherwise,}
\end{array} \right.
\end{equation}
where $x\leftrightarrow y$ is the event that there is an open path 
(i.e., a path of open edges) from $x$ to $y$ in $X(E\setminus\{e\})$.
As a first application of Holley's Inequality, we get from 
(\ref{eq:single_edge_cond_prob}) that
\begin{equation} \label{eq:cond_prob_domination}
\RC^G_{p,q} \bigl(X \in \cdot \bigm| X(E') = \xi\bigr) 
\leqd
\RC^G_{p,q} \bigl(X \in \cdot \bigm| X(E') = \xi'\bigr)
\end{equation}
whenever $E' \subseteq E$ and $\xi \preccurlyeq \xi'$. 

Our next task is to define the random-cluster model on infinite graphs.
Let $G=(V,E)$ be infinite and locally finite. The definition
(\ref{eq:RC_def}) of random-cluster measures
does not work in this case, because there are uncountably
many different configurations $\xi \in \{0,1\}^E$.
Instead, there are two 
other approaches to defining random-cluster measures
on infinite graphs: one via limiting procedures, and the other via local 
specifications, also known as the 
Dobrushin-Lanford-Ruelle (DLR) equations. 
We shall sketch the first approach.

Let $V_1, V_2, \ldots$ be a sequence of finite vertex sets increasing to
$V$ in the sense that $V_1 \subset V_2 \subset \ldots$ and
$\bigcup_{i=1}^\infty V_i = V$. For any finite $K \subseteq V$, define
\[
E(K) := \bigl\{ \ed{ x,y } \in E \st \, x,y \in K\bigr\},
\]
set $E_i := E(V_i)$
and note that $E_1, E_2, \ldots $ increases to $E$ in the same sense that
$V_1, V_2, \ldots$ increases to $V$. Let $\partial V_i$ be the (inner)
boundary of $V_i$ (defined as in (\ref{eq:inner_boundary})). 
Also set $G_i:= (V_i, E_i)$, and let
$\FRC_{p,q}^{G, i}$ be the probability measure on $\{0,1\}^E$ corresponding
to
picking $X\in \{0,1\}^E$ by letting $X(E_i)$ have distribution
$\RC_{p,q}^{G_i}$ and setting $X(e):=0$ for all $e\in E\setminus E_i$.
Since the projection of $\FRC_{p,q}^{G, i}$ on $\{0,1\}^{E\setminus E_i}$
is nonrandom, we can also view $\FRC_{p,q}^{G, i}$
as a measure on $\{0,1\}^{E_i}$, in which case it coincides with
$\RC_{p,q}^{G_i}$. Applying (\ref{eq:cond_prob_domination}) to the graph
$G_i$ with $E':= E_i \setminus E_{i-1}$ and $\xi\equiv 0$ gives
\[
\FRC_{p,q}^{G, i-1} \leqd \FRC_{p,q}^{G, i} \, ,
\]
so that
\begin{equation} \label{eq:monotone_limit}
\FRC_{p,q}^{G, 1} \leqd \FRC_{p,q}^{G, 2} \leqd \cdots \,.
\end{equation}
This implies the existence of a limiting (as $i \rightarrow \infty$)
probability measure $\FRC^G_{p,q}$
on $\{0,1\}^E$. 
This limit is independent of the choice of $\{V_i\}_{i=1}^\infty$, and
we call it the random-cluster measure on $G$ with {\bf free boundary 
condition} (hence the ${\sf F}$ in $\FRC$) and parameters $p$ and $q$. 

Next, define $\WRC_{p,q}^{G, i}$ as the probability measure on $\{0,1\}^E$
corresponding to first setting $X(E \setminus E_i) \equiv 1$, and then
picking $X(E)$ in such a way that
\[
\WRC_{p,q}^{G, i}\bigl(X(E_i) = \xi\bigr) = 
\frac{q^{\kkk(\xi)}}{Z} \prod_{e\in E_i} p^{\xi(e)}(1-p)^{1-\xi(e)}
\]
where $\kkk(\xi)$ is the number of connected components of $\xi$ {\bf that
do not intersect $\partial V_i$}, and $Z$ is again a normalizing constant. 
Similarly as in (\ref{eq:monotone_limit}), we get
\[
\WRC_{p,q}^{G, 1} \geqd \WRC_{p,q}^{G, 2} \geqd \cdots
\]
(with the inequalities reversed compared to (\ref{eq:monotone_limit})), 
and thus also a limiting measure
$\WRC_{p,q}^G$ that we call the random-cluster measure on $G$ with
{\bf wired boundary condition} and parameters $p$ and $q$.

Note that the free and wired random-cluster measures
$\FRC$ and $\WRC$ are both automorphism invariant. This follows from
their construction, in particular from the independence of the choice of 
$\{G_i=(V_i, E_i)\}_{i=1}^\infty$.

\subsection{The Potts model}  \label{subsect:Potts}

Fix a finite graph $G=(V,E)$ and the inverse temperature parameter 
$\beta\geq 0$.
We define the {\bf Gibbs measure for the $q$-state Potts model on $G$
at inverse temperature $\beta$}, denoted $\Pt:=\Pt^G_{q,\beta}$, 
as the probability measure that to each $\omega\in\{1, \ldots, q\}^V$
assigns probability 
\[
\Pt(\omega) :=
\frac{1}{Z} \exp\left( -2\beta \sum_{\ed{ x,y } \in E}
\I{\omega(x) \neq \omega(y)} \right) \, , 
\]
where $Z$ is yet another normalizing constant. 
The main link between random-cluster and Potts
models is the following well-known result.
(See, e.g., \cite{SW}.) 
\begin{prop} \label{prop:from_RC_to_Potts}
Fix a finite graph $G$, an integer $q\geq 2$ and $p\in [0,1]$. Pick
a random edge configuration $X\in \{0,1\}^E$ according to the random-cluster
measure $\RC^G_{p,q}$. Then, for each connected component $\C$ of $X$, 
pick a spin uniformly from $\{1,\ldots, q\}$, and assign this spin to 
all vertices of ${\cal C}$. Do this independently for different
connected components. The $\{1, \ldots, q\}^V$-valued random spin
configuration arising from this procedure is then distributed according
to the Gibbs measure $\Pt^G_{q,\beta}$ for the $q$-state Potts model on $G$
at inverse temperature $\beta:= -\frac{1}{2}\log(1-p)$.  
\end{prop}
This provides the way (mentioned in the introduction)
to reformulate problems about pairwise dependencies
in the Potts model into problems about connectivity probabilities in the
random-cluster model. Aizenman et al.\ \cite{ACCN} were the first to
exploit such ideas
to obtain results about the phase transition behavior of the Potts model,
and the technique has been of much use since then. 

The case of infinite graphs is slightly more intricate. 
Let $G=(V,E)$ be infinite and locally finite, and
let $\{G_i:=(V_i, E_i)\}_{i=1}^\infty$ be as 
in Section \ref{subsect:random-cluster}. 
For $q\in \{2,3,\ldots\}$ and $\beta \geq 0$, 
define probability measures $\left\{\FPt_{q,
\beta}^{G,i}\right\}_{i=1}^\infty$
on $\{1,\ldots, q\}^V$ in such a way that the projection of
$\FPt_{q, \beta}^{G,i}$ on $\{1,\ldots, q\}^{V_i}$ equals
$\Pt_{q,\beta}^{G_i}$, and the spins on
$V\setminus V_i$ are i.i.d.\ uniformly distributed on $\{1, \ldots, q\}$
and independent of the spins on $V_i$. 
Using Proposition \ref{prop:from_RC_to_Potts}, one 
can show that $\FPt_{q, \beta}^{G,i}$
has a limiting distribution $\FPt_{q, \beta}^G$ as $i \rightarrow \infty$. 

Furthermore, for a fixed spin $r\in \{1, \ldots, q\}$, define 
$\WPt_{q, \beta, r}^{G,i}$ to be the distribution corresponding to
picking $X\in \{1, \ldots, q\}^{V}$ by letting $X(V\setminus V_i) \equiv r$,
and letting $X(V_i)$ be distributed according to 
$\Pt_{q,\beta}^{G_i}$ {\bf conditioned on the event that
$X(\partial V_i)\equiv r$}. Again, it turns out that
$\WPt_{q, \beta, r}^{G,i}$ has a limiting
distribution as $i\rightarrow \infty$, and we denote it by
$\WPt_{q, \beta,r}^G$. 

The existence of the limiting distributions
$\FPt_{q, \beta}^G$ and $\WPt_{q, \beta,r}^G$ are nontrivial results,
and in fact the shortest route to proving them goes via random-cluster
arguments: First carry out
the stochastic 
monotonicity arguments for the random-cluster model outlined in Section
\ref{subsect:random-cluster}, and then use Propositions 
\ref{prop:from_FRC_to_FPotts} and \ref{prop:from_WRC_to_WPotts} below. 

A probability measure $\mu$
on $\{1, \ldots, q\}^V$ is said to be a Gibbs measure
(in the DLR sense) for the $q$-state Potts model on $G$ at inverse 
temperature $\beta$, if it admits conditional distributions such that
for all $v\in V$, all $r\in \{1, \ldots, q\}$, and all 
$\omega\in \{1, \ldots, q\}^{V\setminus \{v\}}$, we have 
\rlabel DLR_Potts
{\mu\bigl(X(v)=r \bigm| X(V\setminus \{v\}) = \omega\bigr) = \frac{1}{Z}
\exp\Bigl( -2\beta \sum_{\ed{ v,y } \in E}
\I{\omega(y) \neq r} \Bigr)\,,
}
where the normalizing constant $Z$ may depend on $v$ and $\omega$ but not 
on $r$. The limiting measures $\FPt_{q, \beta}^G$ and $\WPt_{q, \beta, r}^G$ 
are both Gibbs measures in this sense. 

The following extensions of Proposition \ref{prop:from_RC_to_Potts}
provide the relations between $\FRC$ and $\WRC$ on one hand, and
$\FPt$ and $\WPt$ on the other. 
\begin{prop} \label{prop:from_FRC_to_FPotts}
Let $G$ be an infinite locally finite graph, and fix 
$q\in \{2,3,\ldots\}$ and 
$p\in [0,1]$. Pick
a random edge configuration $X\in \{0,1\}^E$ according to $\FRC^G_{p,q}$. 
Then, for each connected component $\C$ of $X$ independently, 
pick a spin uniformly from $\{1,\ldots, q\}$, and assign this spin to 
all vertices of ${\cal C}$. The $\{1, \ldots, q\}^V$-valued random spin
configuration arising from this procedure is then distributed according
to the Gibbs measure $\FPt^G_{q,\beta}$ for the $q$-state Potts model on $G$
at inverse temperature $\beta:= -\frac{1}{2}\log(1-p)$.  
\end{prop}
\begin{prop} \label{prop:from_WRC_to_WPotts}
Let $G$, $p$ and $q$ be as in Proposition \ref{prop:from_FRC_to_FPotts}.
Pick
a random edge configuration $X\in \{0,1\}^E$ according to the random-cluster
measure $\WRC^G_{p,q}$. Then, for each {\bf finite} 
connected component $\C$ of $X$ independently, 
pick a spin uniformly from $\{1,\ldots, q\}$, and assign this spin to 
all vertices of ${\cal C}$. Finally assign value $r$ to all vertices of
infinite connected components. The $\{1, \ldots, q\}^V$-valued random spin
configuration arising from this procedure is then distributed according
to the Gibbs measure $\WPt^G_{q,\beta, r}$ for the $q$-state Potts model on
$G$ at inverse temperature $\beta:= -\frac{1}{2}\log(1-p)$.  
\end{prop}

\section{A dynamical construction} \label{sect:dynamics}

Let $G=(V,E)$ be infinite and locally finite, and let 
$\{G_i:=(V_i, E_i)\}_{i=1}^\infty$ be as in Section \ref{sect:preliminaries}. 
We know from Section \ref{subsect:random-cluster} that
\begin{equation} \label{eq:stoch_dom_1}
\FRC_{p,q}^{G,1} \leqd \FRC_{p,q}^{G,2} \leqd \cdots \leqd
\FRC_{p,q}^G \leqd \WRC_{p,q}^G \leqd \cdots \leqd  
\WRC_{p,q}^{G,2} \leqd \WRC_{p,q}^{G,1}
\,.
\end{equation}
Other well-known stochastic inequalities are that for $p_1 \leq p_2$ and
$i\in \{1,2,\ldots\}$, we have
\begin{equation} \label{eq:stoch_dom_2}
\FRC_{p_1,q}^{G,i} \leqd \FRC_{p_2,q}^{G,i} \, ,
\end{equation}
\begin{equation} \label{eq:stoch_dom_3}
\FRC_{p_1,q}^G \leqd \FRC_{p_2,q}^G \, ,
\end{equation}
\begin{equation} \label{eq:stoch_dom_4}
\WRC_{p_1,q}^{G,i} \leqd \WRC_{p_2,q}^{G,i} \, ,
\end{equation}
and
\begin{equation} \label{eq:stoch_dom_5}
\WRC_{p_1,q}^G \leqd \WRC_{p_2,q}^G \, .
\end{equation}
For all of the above stochastic inequalities, it is desirable to find some 
natural construction of couplings that witness them. What we shall construct
in this section is a coupling of all of the above probability measures
(for all $p\in [0,1]$, $q\geq 1$ and $i\in \{1,2,\ldots \}$) 
{\em simultaneously} that provides witnesses to the 
stochastic inequalities (\ref{eq:stoch_dom_1})--(\ref{eq:stoch_dom_5})
above. Some additional useful aspects of the construction are the
following.
\begin{itemize}
\item[(A1)]
Not only are $\FRC$ and $\WRC$ automorphism 
invariant separately, but also their joint behavior 
in our coupling is automorphism invariant. This remains 
true also if we consider the realizations simultaneously for different
parameter values. See Section \ref{subsect:FRC_crit}, where we describe
an application
where this property is crucial. 
\item[(A2)]
If $G$ is obtained as an automorphism-invariant percolation process on
another
graph $H$, then the construction is easily set up in such a way that
the joint distribution of $G$ and the random-cluster measures on $G$ becomes
an automorphism-invariant process on $H$. (See \cite{HSS} for an example
where an analogous property turns out to be important in the context of
Ising models with external field on percolation clusters.) 
\end{itemize}
Nevertheless, there are still some desirable aspects of couplings of 
random-cluster processes for which we do not know whether or not they hold
for
our construction; see Conjecture \ref{conj:simul_uniqueness} and
Question \ref{quest:cond_prob_gap} in the final section. 

The construction is based on time dynamics for the random-cluster model.
Such time dynamics have previously been considered, e.g., by
Bezuidenhout, Grimmett and Kesten \cite{BGK} and by 
Grimmett \cite{Gr} for the random-cluster
model on $\Z^d$. To some extent our construction will resemble
Grimmett's analysis. However, one feature of our construction that differs
from Grimmett's is that the dynamics are run ``from the past'' rather than
``into the future'', along the lines of the very fashionable
CFTP (coupling from the past) algorithm of Propp and Wilson \cite{PW};
see also \cite{Th} for an early treatment of dynamics from the
past, and \cite{DF} for a survey putting the ideas in a more general
mathematical context.
For the case of finite graphs,  CFTP was applied to simulate
the random-cluster model
in \cite{PW}. Simulation on infinite graphs would require additional
arguments, but our purpose is not simulation; rather, it is to gain some
theoretical information.  For
models other than the random-cluster model, 
CFTP ideas have been extended to the setting of infinite graphs
in \cite{vdBS}, \cite{HS} and \cite{HSS}, but in all those cases the
interaction of the dynamics had a strictly local character, which is not
the case in our context. Another feature of our construction is the
simultaneity in the parameter space. Such simultaneity, which is related to
the level-set representations of Higuchi \cite{Hi}, appears in both
\cite{Gr} and \cite{PW}; Propp and Wilson use the term ``omnithermal''
to denote this particular feature of the construction. 

Let us start with a simple finite case: how do we construct a 
$\{0,1\}^{E_i}$-valued random element with distribution $\RC^{G_i}_{p,q}$
(equivalently, with distribution $\FRC^{G, i}_{p,q}$)? If we are content
with getting something that has only approximately the right distribution,
then the following dynamical approach works fine: Define some
ergodic Markov chain whose unique equilibrium distribution is
$\RC^{G_i}_{p,q}$, and run it for time $T$ starting from
an arbitrary initial state $\xi$. If $T$ is large enough, then the 
distribution of the final state is close to $\RC^{G_i}_{p,q}$, regardless
of the choice of $\xi$. 

In particular, we may proceed as follows.
To each edge $e \in E_i$, we independently
assign an i.i.d.\ sequence $(\phi^e_1, \phi^e_2, \ldots)$ of exponential
random
variables with mean $1$, and an independent i.i.d.\  
sequence $(U^e_1, U^e_2, \ldots)$ of uniform $[0,1]$ random variables. 
For $e\in E_i$ and $k=1,2,\ldots$, let 
$\tau^e_k:= \phi^e_1+ \ldots + \phi^e_k$,
so that $(\tau^e_1, \tau^e_2, \ldots)$ are the jump times of a unit rate
Poisson process. Now define a $\{0,1\}^{E_i}$-valued continuous-time
Markov chain $\{^\xi \bar{X}_{p,q}^{G_i}(t)\}_{t\geq 0}$ with starting state
$^\xi \bar{X}_{p,q}^{G_i}(0):=\xi$ and evolution as follows.
For $e := \ed{ x,y } \in E_i$, 
the value of $^\xi \bar{X}_{p,q}^{G_i}(t)(e)$ does not change
other than (possibly) at the times $\tau^e_1, \tau^e_2, \ldots$, at which
times it takes the value
\begin{equation} \label{eq:evolution_to_the_future}
^\xi \bar{X}_{p,q}^{G_i}(\tau^e_k)(e) :=
\left\{
\begin{array}{ll}
1 & \mbox {if }U^e_k < p \mbox{ and $x\leftrightarrow y$ in } 
{^\xi \bar{X}_{p,q}^{G_i}}(\tau^e_k)(E_i\setminus \{e\}) \\
1 &\phantom{\Bigm|} \mbox{if }U^e_k < \frac{p}{p+(1-p)q} \mbox { and } \neg
\left(
x\leftrightarrow y
\mbox{ in } {^\xi \bar{X}_{p,q}^{G_i}}(\tau^e_k)(E_i\setminus \{e\}) \right)
\\
0 & \mbox{otherwise,}
\end{array} \right. 
\end{equation}
where $\neg$ denotes negation.
(Note that a.s., $\tau^e_k \ne \tau^{e'}_j$ for all $j, k$ when $e \ne
e'$.)
It is easy to see that
this Markov chain is irreducible and reversible with $\RC^{G_i}_{p,q}$ as
stationary distribution, so that indeed $^\xi \bar{X}_{p,q}^{G_i}(t)$
converges
in distribution to $\RC^{G_i}_{p,q}$ as $t\rightarrow\infty$. 
Note also that since $p \geq \frac{p}{p+(1-p)q}$,
the chain preserves the partial order
$\preccurlyeq $ on $\{0,1\}^{E_i}$; in other words, for all $t\geq 0$ we
have
\begin{equation} \label{eq:preserves_partial_order}
^\xi \bar{X}_{p,q}^{G_i}(t) \preccurlyeq {^\eta \bar{X}_{p,q}^{G_i}}(t)
\; \mbox{ whenever } \; \xi \preccurlyeq \eta \, .
\end{equation}
To get a $\{0,1\}^{E_i}$-valued random object whose distribution is 
precisely $\RC^{G_i}_{p,q}$, we need to consider some limit as 
$t \rightarrow \infty$. On the other hand, $^\xi \bar{X}_{p,q}^{G_i}(t)$
does 
not converge in any a.s.\ sense, so this may appear not to be feasible. 

The solution, which turns the convergence in distribution into a.s.\ 
convergence, is to run the dynamics from the past up to time $0$, rather
than from time $0$ into the future. For $T\geq 0$, define the 
$\{0,1\}^{E_i}$-valued continuous-time Markov chain 
\[
\left\{ {^\f_{\; -T}X^{G_i}_{p,q}(t)} \right\}_{t\in [-T, 0]}
\]
with starting state $^\f_{\; -T}X^{G_i}_{p,q}(-T) \equiv 0$ and the
following
evolution, similar to the one of $^\xi \bar{X}_{p,q}^{G_i}$. The value at
an edge $e := \ed{ x,y } \in E_i$ changes only at times
$( \ldots, -\tau_2^e, -\tau_1^e)$, when it takes the value
\begin{equation}  \label{eq:free_dynamics}
_{\; -T}^\f X_{p,q}^{G_i}(-\tau^e_k)(e) :=
\left\{
\begin{array}{ll}
1 & \mbox {if }U^e_k < p \mbox{ and $x\leftrightarrow y$ in } 
{_{\; -T}^\f X_{p,q}^{G_i}}(-\tau^e_k)(E_i\setminus \{e\}) \\
1 & \phantom{\bigg|}\mbox{if }U^e_k < \frac{p}{p+(1-p)q} \mbox { and } \neg
\left(
x\leftrightarrow y
\mbox{ in } {_{\; -T}^\f X_{p,q}^{G_i}}(-\tau^e_k)(E_i\setminus \{e\}) 
\right) \\
0 & \mbox{otherwise,}
\end{array} \right. 
\end{equation}
as in (\ref{eq:evolution_to_the_future}). 
We have, for $0 \leq T_1 \leq T_2$, that
\[
^\f_{\; -T_1}X^{G_i}_{p,q}(0) \preccurlyeq {^\f_{\; -T_2}X^{G_i}_{p,q}}(0)
\]
(essentially because of (\ref{eq:preserves_partial_order})), so by 
monotonicity ${^\f_{\; -T}X^{G_i}_{p,q}}(0)$ has an a.s.\ limit
$^\f X^{G_i}_{p,q}\in \{0,1\}^{E_i}$, defined by setting
$\phantom{\bigg|}^\f X^{G_i}_{p,q}(e) := \lim_{T \rightarrow \infty}
{^\f_{\; -T}X^{G_i}_{p,q}}(0)(e)$ for each $e \in E_i$. Clearly, 
${^\f_{\; -T}X^{G_i}_{p,q}}(0)$ has the same distribution as
$^\xi \bar{X}_{p,q}^{G_i}(T)$ with $\xi \equiv 0$, so
${^\f_{\; -T}X^{G_i}_{p,q}}(0)$ converges in distribution to
$\RC^{G_i}_{p,q}$
as $T\rightarrow \infty$. Hence $\phantom{\bigg|}^\f X^{G_i}_{p,q}$ has
distribution
$\RC^{G_i}_{p,q}$, and if we furthermore define
$^\f X^{G,i}_{p,q}\in \{0,1\}^E$ by setting
\[
^\f X^{G,i}_{p,q}(e) := \left\{
\begin{array}{ll} 
^\f X^{G_i}_{p,q}(e) & \mbox{for } e\in E_i \\
0 & \mbox{otherwise}
\end{array} \right. 
\]
for each $e\in E$, then $^\f X^{G,i}_{p,q}$ has distribution
$\FRC_{p,q}^{G,i}$. 

Now suppose that we have defined the random variables 
$(\phi^e_1, \phi^e_2, \ldots)$ and $(U^e_1, U^e_2, \ldots)$ 
for all $e\in E$ (and not just all $e\in E_i$) in the obvious way. 
By another application of the order-preserving property 
(\ref{eq:preserves_partial_order}), we get that
\[
^\f X^{G,1}_{p,q} \preccurlyeq {^\f X^{G,2}_{p,q}} \preccurlyeq \ldots
\]
so that the limiting object $^\f X^G_{p,q}$, defined by taking
$^\f X^G_{p,q}(e) := \lim_{i\rightarrow\infty} {^\f X^{G,i}_{p,q}} (e)$,
exists. 
For any cylinder set $A \in \{0,1\}^E$, we have
\begin{equation} \label{eq:cylinders}
\P\left({^\f X^G_{p,q}}\in A\right) = 
\lim_{i \rightarrow \infty}\P\left({^\f X^{G,i}_{p,q}}\in A\right) =
\lim_{i \rightarrow \infty}\FRC^{G, i}_{p,q}(A) =
\FRC^G_{p,q}(A)
\end{equation}
so that ${^\f X^G_{p,q}}$ has distribution $\FRC^G_{p,q}$. 
Thus, to summarize the construction so far, what we have is a coupling of
$\{0,1\}^E$-valued random objects 
$^\f X_{p,q}^{G,1}, {^\f X_{p,q}^{G,2}}, \ldots$ and
$^\f X_{p,q}^G$ that witnesses the stochastic inequalities in the
first half of (\ref{eq:stoch_dom_1}). 

Next, we go on to construct, in analogous fashion, 
the corresponding objects for wired
random-cluster measures. For $T\geq 0$, define the $\{0,1\}^E$-valued
continuous-time Markov chain
\[
\left\{ {^\w_{\;\;\;-T} X^{G, i}_{p,q} (t)} \right\}_{t\in [-T,0]}
\] 
with starting configuration $^\w_{\;\;\;-T}X_{p,q}^{G,i}(-T) \equiv 1$. 
Edges $e\in E\setminus E_i$ remain in state $1$ forever, while the value
of an edge $e:= \ed{ x,y } \in E_i$ is updated at times 
$(\ldots, -\tau_2^e, -\tau_1^e)$, when it takes the value
\begin{equation} \label{eq:wired_dynamics}
_{\;\;\; -T}^\w X_{p,q}^{G,i}(-\tau^e_k)(e) :=
\left\{
\begin{array}{ll}
1 & \mbox {if }U^e_k < p 
\mbox{ and } A(x, y, i, p, q, e, k) \\
1 & \mbox{if }U^e_k < \frac{p}{p+(1-p)q} \mbox { and } \neg 
A(x, y, i, p, q, e, k) \\
0 & \mbox{otherwise;}
\end{array} \right. 
\end{equation}
here, $A(x, y, i, p, q, e, k)$ is the event
$\left\{x\stackrel{\partial V_i}{\longleftrightarrow} y \mbox{ in }
{_{\;\;\; -T}^\w X_{p,q}^{G,i}}(-\tau^e_k)(E_i\setminus \{e\})\right\}$,
where, in turn, $x\stackrel{\partial V_i}{\longleftrightarrow} y$ 
denotes the event that
either
\begin{description}
\item{(a) } there is an open path from $x$ to $y$ (not using $e$), or
\item{(b) } both $x$ and $y$ have open paths (not using $e$) to
$\partial V_i$. 
\end{description}
It is immediate from the definition of $\WRC_{p,q}^{G,i}$ that the
conditional $\WRC_{p,q}^{G,i}$-probability that an edge 
$e:= \ed{ x,y } \in E_i$ is open, given the status of all other edges,
is $p$ or ${p}/[{p+(1-p)q}]$, depending on whether or not the event
$\phantom{\bigm|}x\stackrel{\partial V_i}{\longleftrightarrow} y$ happens.
It follows
that the distribution of $^\w_{\;\;\;-T}X_{p,q}^{G,i}(0)$ tends to
$\WRC_{p,q}^{G,i}$ as $T \rightarrow \infty$. Moreover, the dynamics in
(\ref{eq:wired_dynamics}) preserves $\preccurlyeq$ similarly as in
(\ref{eq:preserves_partial_order}), implying that
\[
{_{\;\;\; -T_1}^\w X_{p,q}^{G,i}} \succcurlyeq {_{\;\;\; -T_2}^\w
X_{p,q}^{G,i}}
\]
whenever $0 \leq T_1 \leq T_2$. This establishes the existence of a limiting
$\{0,1\}^E$-valued random object $\phantom{\Big|}^\w X_{p,q}^{G,i}$ defined
by
$^\w X_{p,q}^{G,i}(e) := 
\lim_{T\rightarrow \infty} {^\w X_{p,q}^{G,i}}(0)(e)$ for
each $e\in E$. Clearly, $\phantom{\Big|}^\w X_{p,q}^{G,i}$ has distribution
$\WRC_{p,q}^{G,i}$. Another use of the $\preccurlyeq$-preserving property of
the dynamics (\ref{eq:wired_dynamics}) shows that
\[
^\w X_{p,q}^{G,1} \succcurlyeq {^\w X_{p,q}^{G,2}}
\succcurlyeq \cdots
\,,
\]
so that we have a limiting object $^\w X_{p,q}^G \in \{0,1\}^E$
defined by setting $^\w X_{p,q}^G(e):= \lim_{i \rightarrow \infty}
{^\w X_{p,q}^{G,i}}(e)$ for each $e\in E$. By arguing as in 
(\ref{eq:cylinders}), we get that $^\w X_{p,q}^G$ has
distribution $\WRC_{p,q}^G$. The random objects
$^\w X_{p,q}^{G,1}, {^\w X_{p,q}^{G,2}} ,\ldots$ and
$^\w X_{p,q}^G$ witness the stochastic inequalities in the
second half of (\ref{eq:stoch_dom_1}). 

In order to fully establish that we have a witness to
(\ref{eq:stoch_dom_1}), 
it remains to show that $^\f X_{p,q}^G$ and 
$^\w X_{p,q}^G$ witness the middle inequality in 
(\ref{eq:stoch_dom_1}), i.e., we need to show that 
$\phantom{\Big|}^\f X_{p,q}^G \preccurlyeq {^\w X_{p,q}^G}$. From the
observations that the right-hand sides of (\ref{eq:free_dynamics}) and
(\ref{eq:wired_dynamics}) are increasing in the configurations on 
$E_i \setminus \{e\}$, and that for each such configuration the right-hand 
side of (\ref{eq:wired_dynamics}) is greater than that of
(\ref{eq:free_dynamics}), we get that
\[
^\f_{\; -T}X_{p,q}^{G,i} (t) \preccurlyeq {^\w_{\;\;\; -T}X_{p,q}^{G,i}} (t)
\]
for any $i \in \{1,2,\ldots \}$, $T\geq 0$ and $t\in [-T, 0]$. By taking
$t:=0$,
letting $T\rightarrow \infty$ and then $i \rightarrow \infty$, we get
\begin{equation} \label{eq:witness_1}
^\f X_{p,q}^G \preccurlyeq {^\w X_{p,q}^G} 
\end{equation}
as desired. Hence our coupling is a witness to all the inequalities in
(\ref{eq:stoch_dom_1}). 

It remains to be demonstrated that the coupling is also a witness to the 
inequalities (\ref{eq:stoch_dom_2})--(\ref{eq:stoch_dom_5}). Note first that
the right-hand sides of (\ref{eq:free_dynamics}) and
(\ref{eq:wired_dynamics}) are increasing not only in the configurations on 
$E_i \setminus \{e\}$, but also in $p$. It follows that for $p_1 \leq p_2$
we have
\[
^\f_{\; -T}X_{p_1,q}^{G,i} (t) \preccurlyeq {^\f_{\; -T}X_{p_2,q}^{G,i}} (t)
\]
and
\[
^\w_{\;\;\; -T}X_{p_1,q}^{G,i} (t) \preccurlyeq 
{^\w_{\;\;\; -T}X_{p_2,q}^{G,i}} (t) 
\]
for all $i\in \{1,2,\ldots\}$, $T\geq 0$ and $t\in [-T, 0]$. Taking $t:=0$
and
letting $T \rightarrow \infty$ yields 
\[
^\f X_{p_1,q}^{G,i}  \preccurlyeq  {^\f X_{p_2,q}^{G,i}}  
\]
and
\[
^\w X_{p_1,q}^{G,i}  \preccurlyeq  {^\w X_{p_2,q}^{G,i}} \, ,  
\]
witnessing (\ref{eq:stoch_dom_2}) and (\ref{eq:stoch_dom_4}). Letting
$i \rightarrow \infty$, we get 
\begin{equation} \label{eq:witness_3}
^\f X_{p_1,q}^G  \preccurlyeq  {^\f X_{p_2,q}^G}  
\end{equation}
and
\begin{equation} \label{eq:witness_5}
^\w X_{p_1,q}^G  \preccurlyeq  {^\w X_{p_2,q}^G} \, ,  
\end{equation}
finally witnessing (\ref{eq:stoch_dom_3}) and (\ref{eq:stoch_dom_5}).
In fact, examination also shows that as long as $p_1 \le p_2$ and
$p_1/[(1-p_1)q_1] \le p_2/[(1-p_2)q_2]$, we have
\begin{equation} \label{eq:witness_genfree}
^\f X_{p_1,q_1}^G  \preccurlyeq  {^\f X_{p_2,q_2}^G}  
\,,
\end{equation}
\begin{equation} \label{eq:witness_genwired}
^\w X_{p_1,q_1}^G  \preccurlyeq  {^\w X_{p_2,q_2}^G} 
\,,
\end{equation}
and
\begin{equation} \label{eq:witness_genfreewired}
^\f X_{p_1,q_1}^G  \preccurlyeq  {^\w X_{p_2,q_2}^G} \, ,  
\end{equation}
witnessing more general well-known stochastic inequalities \cite{F}.

Property (A1) of the coupling is obvious from the construction.
In order for (A2) to be true, we need only define random variables
$\{ \phi_k^e, U_k^e\}_{e\in E(H), i=1,2,\ldots}$ for all edges in $H$ and
to take them to
be independent of the percolation process that yields $G$ from $H$.

\section{Bernoullicity} \label{sect:Bernoullicity}

Let $\Gamma$ be a closed subgroup of $\Aut(G)$ with $G = (V, E)$ being any
connected graph.
We shall be most interested in two cases: (1) that $G$ is the Cayley graph of
$\Gamma$ with respect to some finite generating set of
$\Gamma$; and (2) that $\Gamma = \Aut(G)$ and $G$ is
quasi-transitive.
Let
$S$ and $T$ be arbitrary state spaces. For $\gamma\in \Gamma$, define the 
map $\theta_\gamma: S^V \rightarrow S^V$ 
(or $\theta_\gamma: T^V \rightarrow T^V$) by setting
$\theta_\gamma\omega(x) := \omega(\gamma^{-1}x)$ for each $x\in V$.
A measurable mapping $f: (S^V, \mu) \rightarrow (T^V, \nu)$ is said to be
{\bf
$\Gamma$-equivariant} if it commutes with these actions of $\Gamma$, i.e.,
if $f(\theta_\gamma\omega) = \theta_\gamma\bigl(f(\omega)\bigr)$ for all
$\gamma\in \Gamma$ and $\mu$-a.e.\ $\omega \in S^V$; it is called {\bf
measure-preserving} if $\nu = \mu \circ f^{-1}$. The action of $\Gamma$ on
$(T^V, \nu)$ is called {\bf free} if for $\nu$-a.e.\ $x \in T^V$, the only
element in $\Gamma$ that leaves $x$ fixed is the identity.

We say that a probability measure $\nu$ on $T^V$ is a {\bf $\Gamma$-factor
of
an i.i.d.\ process} if there exists a $T^V$-valued random element $X$
with distribution $\nu$, a state space $S$, an $S^V$-valued random element
$Y$ with distribution $\mu$, and a $\Gamma$-equivariant measure-preserving
mapping $f: (S^V, \mu) \rightarrow (T^V, \nu)$ such that
\begin{itemize}
\item[(i)] $Y$ is an i.i.d.\ process, and
\item[(ii)] $X=f(Y)$. 
\end{itemize}
In case $G$ is the Cayley graph of $\Gamma$, if
$S$ can be taken to be finite and
$f$ can be taken to be an invertible
mapping, then $(\Gamma, \nu)$ is said to be {\bf Bernoulli},
a mixing property of
fundamental importance in ergodic theory. 
In \cite[p.~127]{OW}, it is shown that the following definition is a proper
extension of the preceding definition:
An action $(\Gamma, \nu)$ is said to be {\bf Bernoulli}
if it is a free $\Gamma$-factor of a Poisson process on $\Gamma$.
We shall prove, using the
dynamical
construction in Section \ref{sect:dynamics}, that Bernoullicity
holds for the wired Potts model on $\Z^d$, and more generally on many
amenable quasi-transitive graphs.
We shall need the following condition.
Let $S_n(x)$ denote the set of points at distance $n$ from a vertex $x$.
Consider the condition on $\Gamma$ that
\rlabel tobefree
{\forall x \in V \quad \forall y \in \Gamma x \setminus \{ x \}
\quad \exists \mbox{ infinitely many } n \qquad 
S_n(x) \ne S_n(y)
\,.
}
\begin{thm} \label{thm:Bernoullicity}
Let $G$ be a Cayley graph of any amenable group $\Gamma$
or be any amenable graph with a closed automorphism group $\Gamma$ acting
quasi-transitively on $G$ and satisfying \rref e.tobefree/.
Let $q\in \{2,3,\ldots\}$, $r\in \{1,\ldots, q\}$, and $\beta\geq 0$.
Then the Gibbs measure 
$\WPt^G_{q, \beta, r}$ is Bernoulli with respect to the action of $\Gamma$. 
\end{thm}

For the $\Z^d$ case, this was previously known only for the cases where
either
$q=2$ (the Ising model) or $\beta$ is sufficiently small; see, e.g., 
\cite{OW2}, \cite{LGR} and \cite{St}. For the Ising model
result on amenable graphs, see \cite{Ad}, while for a proof of a stronger
property than Bernoullicity in the case of $\beta$ small, using CFTP ideas, 
see \cite{HS}. 
The paper \cite{HSS} uses ideas similar to ours to prove that 
the Ising model is Bernoulli.

\procl r.freeness
Actually, we shall prove a slightly stronger result, which is the best
possible. That is, we shall show that as
long as i.i.d.\ variables on the vertices of $G$ yield a free action of
$\Gamma$, then $\WPt^G_{q, \beta, r}$ is Bernoulli.
It is not clear when the full automorphism group $\Aut(G)$ satisfies this
freeness condition, so we have supplied the condition \rref e.tobefree/.
\endprocl

We call an i.i.d.\ process $(S^V, \mu)$ {\bf standard} if $S$ is a standard
Borel space and the marginal of $\mu$ on $S$ is Borel.
Ornstein and Weiss \cite{OW} show that when $\Gamma$ is amenable and
discrete, then $(\Gamma, \nu)$ is Bernoulli iff it is a free
$\Gamma$-factor of a standard i.i.d.\ process.
More generally, we have the following result:

\procl l.AutBern
Let $V$ be a countable set and $\Gamma$ be a closed subgroup of the
symmetric group on $V$.
Suppose that all orbits of the $\Gamma$-action on $V$ are infinite
and that $\Gamma$ is amenable, unimodular, and not the union of an
increasing sequence of compact proper subgroups of $\Gamma$.
Further, suppose that for each $x\in V$, the $\Gamma$-stabilizer of $x$ is
compact.
Then every free $\Gamma$-factor of a standard i.i.d.\ process $(S^V, \mu)$
is Bernoulli.
\endprocl

\proof
Assume that there is some free $\Gamma$-factor $\nu$ of a standard i.i.d.\
process $(S^V, \mu)$, since otherwise there is nothing to prove. 
Let $Z_n$ be i.i.d.\ Poisson point processes on $\Gamma$
with Haar measure as the underlying intensity measure. 
By \cite[Theorem III.6.5]{OW}, the product process $\Seq{Z_n \st n \ge 1 }$
is Bernoulli.
We shall show that $\nu$ is a $\Gamma$-factor of $\Seq{Z_n \st n \ge 1 }$,
whence is a factor of a Poisson process, whence is Bernoulli.

Let $W$ be a selection of one point from each orbit of the action of
$\Gamma$ on $V$.
Given $v \in V$, let $X_n(v)$ be the number of points in $Z_n$ that take
$o$ to $v$ for $v \in V$, where $\{o\} = W \cap \Gamma v$.
Since $\Gamma$ is a countable union of translates of stabilizers, each
stabilizer has positive finite Haar measure, so that
$X_n(v)$ is a nontrivial Poisson random variable.
Also, the random variables $\Seq{X_n(v) \st n \ge 1, v \in V}$ are mutually
independent.
Since $X_n$ is a $\Gamma$-factor of $Z_n$, it follows that
$\Seq{X_n\st n \ge 1}$ is a $\Gamma$-factor of $\Seq{Z_n}$.
Since every standard i.i.d.\ process $(S^V, \mu)$ is a $\Gamma$-factor of
$\Seq{X_n \st n \ge 1}$ and $\nu$ is a factor of $(S^V, \mu)$, we obtain
the result we want.
\Qed

We also need the following fact:

\procl l.good-entropy
If $G$ is a quasi-transitive amenable graph, then $\Aut(G)$ is amenable,
unimodular, and not the union of an increasing sequence of compact proper
subgroups.
\endprocl

\proof
$\Aut(G)$ is amenable and unimodular by results of Soardi
and Woess \cite{SoWo} and Salvatori \cite{Sa}; see also \cite{BLPS} for
another proof. Furthermore, in this case $\Aut(G)$ is generated
by, say,
the compact set $\Delta:= \{\gamma \in \Aut(G) \st d(o, \gamma o) \le 2r+1\}$,
where $r$ is such that every vertex of $G$ is within distance $r$ of some
vertex in $\Aut(G) o$ and $d(\cdot, \cdot)$ denotes distance in $G$.
Thus, if $\Gamma_n$ are compact increasing subgroups of
$\Aut(G)$ whose union is $\Aut(G)$, we have $\bigcap_{n \geq 1} 
(\Delta \setminus
\Gamma_n) = \emptyset$, whence for some $n$, we have $\Delta \subseteq
\Gamma_n$. Since $\Delta$ generates $\Aut(G)$, it follows that $\Gamma_n =
\Aut(G)$. 
\Qed

Because of the above, 
Theorem \ref{thm:Bernoullicity} is established once the following
lemma is proved:
\begin{lem} \label{lem:factor}
For any graph $G$, any subgroup $\Gamma$ of $\Aut(G)$,
any $q\in\{2,3,\ldots\}$ and $r\in \{1, \ldots, q\}$,
and any $\beta \geq 0$, the Gibbs measure $\WPt^G_{q, \beta, r}$ is 
a $\Gamma$-factor of a standard i.i.d.\ process.
If either (i) $\Gamma$ is countable and every
element of $\Gamma$ other than the identity moves an infinite number of
vertices or (ii) $\Gamma$ satisfies condition \rref e.tobefree/, 
then the action of $\Gamma$ on $\WPt^G_{q, \beta, r}$ is free.
\end{lem}
\proof
Let the degree of $G$ be $d$. For each $x \in V$, let $N_x = \{Z^x_1,
\ldots,
Z^x_d\}$ be the set of neighbors of $x$ in any fixed order.

Take 
\[
S:= \bigl\{ [0, \infty) \times [0,1] \bigr\}^{\{1,2,\ldots\} \times \{1,
\ldots, d\}}
\times [0, 1]^d \times [0,1] \times \{1, \ldots, q\} \, .
\]
Let
\[
\left\{ \phi^j_k(x), U^j_k(x), U^j_*(x), U^*(x), \sigma(x) 
\st k=1,2,\ldots,\ j=1, \ldots, d,\ x \in V
\right\}
\]
be independent random variables with $\phi^j_k(x)$ exponential of mean 1,
$U^j_k(x)$, $U^j_*(x)$, and $U^{*}(x)$ uniform $[0, 1]$, and
$\sigma(x)$ uniform on $\{1, \ldots, q\}$.
For each $x\in V$, put
\[
Y(x):= \left( \bigl(\phi^j_k(x), U^j_k(x)\bigr)_{k=1,2,\ldots, j=1, \ldots,
d}, \bigl( U^j_*(x) \bigr)_{j = 1, \ldots, d}, U^*(x),
\sigma(x) \right) \, .
\]

Set $p:=1-e^{-2\beta}$, and construct a $\{0,1\}^E$-valued edge
configuration
$X^G_{p,q}$ with distribution $\WRC^G_{p,q}$ by the dynamical construction
in Section \ref{sect:dynamics}, where for each $e\in E$ we take
\begin{equation} \label{eq:from_vertices_to_bonds}
(\phi^e_k, U^e_k)_{k=1,2,\ldots} := \bigl(\phi_k^j(x),
U_k^j(x)\bigr)_{k=1,2,\ldots}\,,
\end{equation}
where $x\in V$ and $j \in \{1, \ldots, d\}$ are chosen in such a way that
$e=\ed{ x, Z^x_j}$, and, if we denote $y := Z^x_j$ and $j'$ is such that $x
= Z^y_{j'}$, then $U^j_*(x) < U^{j'}_*(y)$.
This choice of $x$ and $j$ is a.s.\ unique.

From $X^G_{p,q}$, we obtain the desired spin configuration 
$X\in \{1, \ldots, q\}^V$ with distribution $\WPt^G_{q,\beta, r}$ by
assigning spins to the connected components of $X^G_{p,q}$ as in 
Proposition \ref{prop:from_WRC_to_WPotts}: All vertices in infinite
connected
components in $X_{p,q}^G$ are assigned value $r$, whereas the vertices of
each
finite connected component $\C$ are assigned value $\sigma(x)$, where $x$
is the vertex in $\C$ that minimizes $U^*(x)$. It is obvious that
this mapping $Y \mapsto X$
from $S^V$ to $\{1, \ldots, q\}^V$ is $\Aut(G)$-equivariant, and
that the resulting spin configuration has distribution
$\WPt^G_{q,\beta, r}$. Hence
$\WPt^G_{q,\beta, r}$ is a factor of a standard i.i.d.\ process. 

To see that the action of $\Gamma$ on $\WPt^G_{q, \beta, r}$
is free under the additional hypotheses (i)
stated in the lemma, it suffices to show that for any $\gamma \in \Gamma$
other than the identity, $\P[\theta_\gamma X = X] = 0$.  From the
hypotheses, we may find an infinite set $W$ of vertices such that $\gamma x
\notin W$ for all $x \in W$ and $\gamma x \ne \gamma y$ for distinct $x, y
\in W$.
Because of \rref e.DLR_Potts/, by repeated conditioning we see that
there is some $c < 1$ such that for any
$x_1, \ldots, x_n \in W$, we have $\P[X(x_i) = X(\gamma^{-1} x_i) \mbox{
for all } i = 1, \ldots, n] \le c^n$. Therefore $\P[\theta_\gamma X = X] =
0$.

Consider now the hypothesis (ii). 
Again because of \rref e.DLR_Potts/, there 
is some $c < 1$ such that if $A$ and $A'$ are two finite sets of vertices that
are not identical, then the chance is at most $c$
that the number of spins in $A$ equal to 1
is the same as the number of spins in $A'$ equal to 1, even given all spins
outside $A \cup A'$.
Suppose that $x \ne y$ and $x$ and $y$ are in the same orbit.
Let $W(x, y)$ be the set of spin configurations such that for some $n$,
the number of spins in $S_n(x)$ equal to 1 differs from the number in $S_n(y)$.
By our assumption and the fact just noted, it follows that $W(x, y)$ has
probability 1.
Hence so does $W := \bigcap_{x, y} W(x, y)$.
It is clear that $\Gamma$ acts freely on $W$.
\Qed

\section{Further remarks on the coupling construction}  \label{sect:remarks}

\subsection{Critical behavior of the random-cluster model}
\label{subsect:FRC_crit}

Let us mention another application of the pointwise construction in 
Section \ref{sect:dynamics}.
Consider the random-cluster model on an infinite quasi-transitive graph $G$ at
some fixed value of $q$. We shall let $p$ vary. Clearly, by stochastic 
monotonicity, the $\FRC^G_{p,q}$- and $\WRC^G_{p,q}$-probabilities of having
some infinite open cluster are increasing in $p$. Furthermore, by ergodicity,
these probabilities must be $0$ or $1$ for any given $p$ (although the
$\FRC^G_{p,q}$-probability does not necessarily equal the 
$\WRC^G_{p,q}$-probability). Hence, there exist critical values 
$p_c^\f:=p_c^\f(G,q)$ and $p_c^\w:=p_c^\w(G,q)$ such that
\begin{equation} \label{eq:p_c_free}
\FRC^G_{p,q} (\exists \mbox{ at least one infinite cluster})= \left\{
\begin{array}{ll}
0 & \mbox{for } p<p_c^\f, \\
1 & \mbox{for } p>p_c^\f
\end{array} \right. 
\end{equation}
and
\begin{equation} \label{eq:p_c_wired}
\WRC^G_{p,q} (\exists \mbox{ at least one infinite cluster})= \left\{
\begin{array}{ll}
0 & \mbox{for } p<p_c^\w, \\
1 & \mbox{for } p>p_c^\w \, .
\end{array} \right. 
\end{equation}
A very natural question is whether or not there is an infinite cluster at
criticality. In \cite{HJL}, we proved that when $G$ is a unimodular 
nonamenable quasi-transitive graph,
then the answer is no for $\FRC$. In other words, 
\begin{equation} \label{eq:at_criticality}
\FRC^G_{p_c^\f,q} (\exists \mbox{ at least one infinite cluster})= 0 \, .
\end{equation}
The proof in \cite{HJL} of (\ref{eq:at_criticality}) uses, as a key ingredient,
the existence of an automorphism-invariant coupling of the measures
$\FRC^G_{p,q}$ for different $p$
that witnesses the stochastic domination (\ref{eq:stoch_dom_3}). Such a
coupling was provided in Section \ref{sect:dynamics} of the present paper. 

It seems reasonable to expect that (\ref{eq:at_criticality}) extends to 
all quasi-transitive graphs (except those for which the critical value is $1$). For
$q=1$, this was conjectured by Benjamini and Schramm \cite{BS}. 
The situation for $\WRC$ seems to be more complicated. For instance, 
as shown in \cite{CCST} and \cite{H2}, 
when
$G$ is the regular tree $\T_n$ with $n \geq 2$, we get that the 
$\WRC^G_{p_c^\w,q}$-probability of seeing an infinite cluster is $0$ or $1$
depending on whether $q\in [1,2]$ or $q>2$. 

\subsection{Simultaneity statements}

For quasi-transitive graphs,
the famous finite-energy argument of Newman and Schulman \cite{NS}
shows that the number of infinite clusters must (under either $\FRC$ or
$\WRC$, and for fixed $p$ and $q$) be an almost sure constant, and
either $0$, $1$ or $\infty$. For  
unimodular quasi-transitive graphs, Lyons \cite{L} recently obtained the
necessary uniqueness monotonicity statement for deducing that
(in addition to the critical values in (\ref{eq:p_c_free}) and
(\ref{eq:p_c_wired})), there exist critical values $p_u^\f$ and $p_u^\w$
such that
\begin{equation} \label{eq:p_u_free}
\FRC^G_{p,q} (\exists \mbox{ \rm a unique infinite cluster})= \left\{
\begin{array}{ll}
0 & \mbox{for } p<p_u^\f, \\
1 & \mbox{for } p>p_u^\f
\end{array} \right. 
\end{equation}
and
\begin{equation} \label{eq:p_u_wired}
\WRC^G_{p,q} (\exists \mbox{ \rm a unique infinite cluster})= \left\{
\begin{array}{ll}
0 & \mbox{for } p<p_u^\w, \\
1 & \mbox{for } p>p_u^\w \, .
\end{array} \right. 
\end{equation}
(For $q=1$ this goes back to \cite{HP} and \cite{S}.)
See \cite{HJL} for a detailed discussion of how the four critical values
$p_c^\f$, $p_c^\w$, $p_u^\f$ and $p_u^\w$ relate to each other. 

It is not obvious that, in the coupling of Section \ref{sect:dynamics}, 
(\ref{eq:p_u_free}) and (\ref{eq:p_u_wired}) hold simultaneously 
for all $p$ and $q$. This is in fact an open problem, and
we conjecture the following strengthening, analogous to the simultaneous
uniqueness results of \cite{A}, \cite{HP}, \cite{HPS}, and \cite{S}:

\begin{conj}  \label{conj:simul_uniqueness}
Let $G=(V,E)$ be connected and quasi-transitive.
For a configuration $\xi \in \{0, 1\}^E$, write $N(\xi)$ for the number of
infinite clusters in $\xi$.
Let $\dom$ be the set of quadruples $(p_1, p_2, q_1, q_2)$ such that
$$
p_1 \le p_2 \quad\mbox{ and }\quad
\frac{p_1}{(1-p_1)q_1} \le \frac{p_2}{(1-p_2)q_2}
\,,
$$
with at least one of these inequalities being strict.
In the notation of \rref s.dynamics/, we have a.s.\
for all quadruples $(p_1, p_2, q_1, q_2) \in \dom$ simultaneously,
each infinite cluster of $Y$ contains $N(X)$ infinite clusters of $X$,
where $X$ and $Y$ may be any of the following three pairs
of random variables:
\begin{itemize}
\item[\rm (i) ] 
$X = {^\f X_{p_1,q_1}^G}$ and $Y =  {^\f X_{p_2,q_2}^G}$, 
\item[\rm (ii) ] 
$X = {^\w X_{p_1,q_1}^G}$ and $Y = {^\w X_{p_2,q_2}^G}$,
\item[\rm (iii) ] 
$X = {^\f X_{p_1,q_1}^G}$ and $Y = {^\w X_{p_2,q_2}^G}$.
\end{itemize}
\end{conj}

\subsection{Another open problem}

Let us finally discuss another open problem concerning
our coupling in Section \ref{sect:dynamics}. For $p_1 <p_2$, define
\[
\Delta_q (p_1, p_2) := \min\left\{ p_2-p_1, \,
\frac{p_2}{p_2 +(1-p_2)q} - \frac{p_1}{p_1 +(1-p_1)q} \right\}
\]
and note that $\Delta_q (p_1, p_2) >0$. 
For $e \in E$ and $\xi \in \{0, 1\}^{E \setminus \{e\}}$, write
$A(\xi, e, p, q)$ for the event that
${^\f X_{p,q}^G} (E \setminus \{e\}) = \xi$. 
{}From the fact that
$\FRC^G_{p,q}$ is a DLR random-cluster measure, it follows that for any
$e\in E$ and almost any $(\xi, \eta) \in (\{0,1\}^{E \setminus \{e\}})^2$
with respect to the law of $\left({^\f X_{p_1,q}^G} (E \setminus \{e\}) ,\;
{^\f X_{p_2,q}^G} (E \setminus \{e\})\right)$
under our coupling (which implies that
$\xi \preccurlyeq \eta$), we have
\begin{equation}
\label{eq:cond_prob_difference}
\P\left({^\f X_{p_2,q}^G} (e) = 1 \bigm|
A(\eta, e, p_2, q) \right) -
\P\left({^\f X_{p_1,q}^G} (e) = 1 \bigm| 
A(\xi, e, p_1, q) \right)
\geq \Delta_q (p_1, p_2) 
\end{equation}
(and similarly for wired random-cluster measures; everything we say in 
relation to Question \ref{quest:cond_prob_gap} applies as well
to the wired case as to the free). From this, one is easily seduced
into thinking that 
\begin{equation} \label{eq:coupled_cond_prob_difference}
\P\left( {^\f X_{p_2,q}^G} (e) = 1 \, , \, {^\f X_{p_1,q}^G} (e) = 0 \bigm|
A(\eta, e, p_2, q) \cap A(\xi, e, p_1, q) \right)
\geq \Delta_q (p_1,p_2) \,,
\end{equation}
but to conclude this directly from (\ref{eq:cond_prob_difference}) is
unwarranted, because conditioning on $\xi$ and $\eta$ jointly is not the
same as conditioning on them separately. It is nevertheless natural to ask 
whether something like (\ref{eq:coupled_cond_prob_difference}) is true. In
particular, the following question asks for a weaker property. 
\begin{quest} \label{quest:cond_prob_gap}
For $p_1 < p_2$ and $q \geq 1$, does there exist an $\eps > 0$ 
(depending on $p_1$, $p_2$ and $q$) such that for any $e\in E$ and almost
any $(\xi, \eta)\in (\{0,1\}^{E\setminus\{e\}})^2$, we have
\[
\P\left( {^\f X_{p_2,q}^G} (e) = 1 \, , \, {^\f X_{p_1,q}^G} (e) = 0 \bigm|
A(\eta, e, p_2, q) \cap A(\xi, e, p_1, q) \right)
\geq \eps \, ?
\]
\end{quest}
A positive answer to this question (for our coupling or for
some other automorphism-invariant witness to the stochastic inequality
$\FRC_{p_1,q}^G \leqd \FRC_{p_2,q}^G$) is precisely the missing ingredient
that prevented the authors of \cite{HP} from extending their uniqueness
monotonicity result for i.i.d.\ percolation ($q=1$) for unimodular
quasi-transitive
graphs to the more general case $q\geq 1$ (i.e., from proving 
the relations \rref e.p_u_free/ and \rref e.p_u_wired/ that were later 
obtained in \cite{L}). Such a positive answer
might perhaps also be an ingredient in applying
the reasoning of Schonmann \cite{S} in order to remove the unimodularity 
assumption in these results.

\medbreak
\noindent {\bf Acknowledgement.}\enspace
We are grateful to Benjy Weiss for useful discussions.


\begin{thebibliography}{99}

{\small

\bibitem{Ad} Adams, S. (1992) F\o lner independence and the amenable Ising 
model, {\sl Ergod. Th. Dynam. Sys.} {\bf 12}, 633--657. 

\bibitem{ACCN} Aizenman, M., Chayes, J.T., Chayes, L., and Newman, C.M.
(1988) Discontinuity of the magnetization in one-dimensional
$1/|x-y|^2$ Ising and Potts models, {\sl J. Statist. Phys.} {\bf 50},
1--40.

\bibitem{A} Alexander, K. (1995) Simultaneous uniqueness of infinite
clusters in stationary random labeled graphs, {\sl Commun. Math. Phys.}
{\bf 168}, 39--55.

\bibitem{BLPS} Benjamini, I., Lyons, R., Peres, Y., and Schramm, O. (1999)
Group-invariant percolation on graphs, {\sl Geom. Funct. Analysis} 
{\bf 9}, 29--66.

\bibitem{BS} Benjamini, I. and Schramm, O. (1996) Percolation beyond
${\bf Z}^d$, many questions and a few answers, {\sl Electr. Commun. Probab.}
{\bf 1}, 71--82. 

\bibitem{BS:hp} Benjamini, I. and Schramm, O. (2001) Percolation in
the hyperbolic plane, {\sl J. Amer. Math. Soc.} {\bf 14}, 487--507.  

\bibitem{vdBS} van den Berg, J. and Steif, J.E. (1999) On the existence and
nonexistence of finitary codings for a class of random fields, 
{\sl Ann. Probab.} {\bf 27}, 1501--1522.

\bibitem{BGK} Bezuidenhout, C., Grimmett, G. and Kesten, H. (1993)
Strict inequality for critical values of Potts models and random-cluster 
processes, {\sl Commun. Math. Phys.} {\bf 158}, 1--16.  

\bibitem{BK} Burton, R.M. and Keane, M.S. (1989) Density and uniqueness in
percolation, {\sl Commun. Math. Phys.} {\bf 121}, 501--505.

\bibitem{CCST} Chayes, J.T., Chayes, L., Sethna, J.P., and Thouless, D.J.
(1986).  A mean field spin glass with short-range interactions,  {\sl Comm.
Math. Phys.} {\bf 106}, 41--89.

\bibitem{DF} Diaconis, P. and Freedman, D. (1999) Iterated random functions,
{\sl SIAM Rev.} {\bf 41}, 45--76. 

\bibitem{ES} Edwards, R.G. and Sokal, A.D. (1988) Generalization of the
Fortuin-\-Kasteleyn-\-Swendsen-\-Wang representation and Monte Carlo
algorithm, {\sl Phys. Rev.} D {\bf 38}, 2009--2012.

\bibitem{F} Fortuin, C.M. (1972).  On the random-cluster model. {I}{I}{I}.
{T}he simple random-cluster model,  {\sl Physica} {\bf 59}, 545--570.

\bibitem{FK} Fortuin, C.M. and Kasteleyn, P.W. (1972)
On the random-cluster model. I. Introduction and relation to other models,
{\sl Physica} {\bf 57}, 536--564.

\bibitem{GHM} Georgii, H.-O., H\"aggstr\"om, O. and Maes, C. (2001)
The random geometry of equilibrium phases, {\sl Phase Transitions and
Critical
Phenomena, Volume 14} 
(C. Domb and J.L. Lebowitz, eds), pp 1--142, Academic Press, London.

\bibitem{Gr} Grimmett, G.R. (1995) The stochastic random-cluster process,
and the uniqueness of random-cluster measures, {\sl Ann. Probab.} {\bf 23}, 
1461--1510. 

\bibitem{H2} H\"aggstr\"om, O. (1996) The random-cluster model on a
homogeneous tree, {\sl Probab. Th. Rel. Fields} {\bf 104}, 231--253.

\bibitem{HJL} H\"aggstr\"om, O., Jonasson, J. and Lyons, R. (2001)
Explicit isoperimetric constants and phase transitions in the
random-cluster model, {\sl Ann. Probab.}, to appear. 

\bibitem{HP} H\"aggstr\"om, O. and Peres, Y. (1999) Monotonicity of
uniqueness for percolation on Cayley graphs: all infinite clusters are
born simultaneously, {\sl Probab. Th. Rel. Fields} {\bf 113},
273--285.

\bibitem{HPS} H\"aggstr\"om, O., Peres, Y. and Schonmann, R.  (1999)
Percolation on transitive graphs as a coalescent process: relentless
merging followed by simultaneous uniqueness.  In Bramson, M. and
Durrett, R., editors, {\sl Perplexing Probability Problems: Papers in Honor
of Harry Kesten}, pages 69--90, Boston.  Birkh{\"a}user.

\bibitem{HSS} H\"aggstr\"om, O., Schonmann, R.H. and Steif, J.E. (2000) The
Ising model on diluted graphs and strong amenability, {\sl Ann. Probab.}
{\bf 28}, 1111--1137. 

\bibitem{HS} H\"aggstr\"om, O. and Steif, J.E. (2000) Propp-Wilson
algorithms 
and finitary codings for high noise Markov random fields, 
{\sl Combin. Probab. Computing} {\bf 9}, 425--439. 

\bibitem{Hi} Higuchi, Y. (1991) Level set representation for the Gibbs
states 
of the ferromagnetic Ising model, {\sl Probab. Th. Rel. Fields} {\bf 90},
203--221.

\bibitem{LGR} Liberto, F. di, Gallavotti, G. and Russo, L. (1973)
Markov processes, Bernoulli schemes, and Ising model,
{\sl Commun. Math. Phys.} {\bf 33}, 259--282.

\bibitem{Li} Lindvall, T. (1992) {\sl Lectures on the Coupling Method},
Wiley, New York.

\bibitem{L} Lyons, R. (2000) Phase transitions on nonamenable graphs, 
{\sl J. Math. Phys.} {\bf 41}, 1099--1126.

\bibitem{NS} Newman, C.M. and Schulman, L.S. (1981) Infinite clusters in
percolation models, {\sl J. Statist. Phys.} {\bf 26}, 613--628. 

\bibitem{OW2} Ornstein, D.S. and Weiss, B. (1973) $\Z^d$-actions and
the Ising model, {\sl unpublished manuscript}. 

\bibitem{OW} Ornstein, D.S. and Weiss, B. (1987)
Entropy and isomorphism theorems for actions of amenable
groups, {\sl J. Anal. Math.} {\bf 48}, 1--141.

\bibitem{PW} Propp, J.G. and Wilson, D.B. (1996)
Exact sampling with coupled Markov chains and applications to statistical
mechanics, {\sl Random Structures Algorithms} {\bf 9}, 223--252.

\bibitem{Sa} Salvatori, M. (1992) On the norms of group-invariant
transition operators on graphs, {\sl J. Theoret. Probab.} {\bf 5},
563--576.

\bibitem{S} Schonmann, R.H. (1999) Stability of infinite clusters in
supercritical percolation, 
{\sl Probab.\ Th.\ Rel.\ Fields} {\bf 113}, 287--300.

\bibitem{SoWo} Soardi, P.M. and Woess, W. (1990) Amenability, unimodularity,
and the spectral radius of random walks on infinite graphs, {\sl Math. Z.}
{\bf 205}, 471--486.

\bibitem{St} Steif, J.E. (1991) 
$\overline{d}$-convergence to equilibrium and space-time
              Bernoullicity for spin systems in the $M<\epsilon$ case,
{\sl Ergod. Th. Dynam. Sys.} {\bf 11}, 547--575.

\bibitem{SW} Swendsen, R.H. and Wang, J.-S. (1987)
Nonuniversal critical dynamics in Monte Carlo simulations,
{\sl Phys. Rev. Lett.} {\bf 58}, 86--88.

\bibitem{Th} Thorisson, H. (1988) Backward limits, {\sl Ann. Probab.}
{\bf 16}, 914--924. 

}

\end{thebibliography}
\end{document}